\documentclass[10pt,english,10pt,onecolumn,draftcls]{IEEEtran}
\usepackage[latin9]{inputenc}
\usepackage{geometry}
\usepackage{color}
\usepackage{verbatim}
\usepackage{amsthm}
\usepackage{amsmath}
\usepackage{amssymb}
\usepackage{graphicx}
\usepackage{esint}

\makeatletter
\theoremstyle{plain}
\newtheorem{thm}{\protect\theoremname}
\theoremstyle{plain}
\newtheorem{lem}[thm]{\protect\lemmaname}
\theoremstyle{plain}
\newtheorem{prop}[thm]{\protect\propositionname}
\theoremstyle{plain}
\newtheorem{cor}[thm]{\protect\corollaryname}

\@ifundefined{definecolor}
 {\usepackage{color}}{}
\usepackage{multicol}
\usepackage{flushend}
\usepackage[linesnumbered,ruled,boxed]{algorithm2e}
\usepackage{float}
\usepackage{footnote} 
\IEEEoverridecommandlockouts

\RequirePackage{colortbl, tabularx}
\@ifundefined{comment}{}
  {
   \renewenvironment{comment}
    {
     \par\medskip\noindent
     \tabularx{\columnwidth}{|>{\columncolor[gray]{0.9}}X|}
     \hline
     \emph{\textbf{Comment:}}
    }
    {
     \endtabularx\hrule\par\medskip
    }
  }%

\makeatother

\usepackage{babel}
\providecommand{\corollaryname}{Corollary}
\providecommand{\lemmaname}{Lemma}
\providecommand{\propositionname}{Proposition}
\providecommand{\theoremname}{Theorem}

\begin{document}

\title{\vspace{16pt}
}

\title{Online Distributed Optimization \hspace{6cm}on Dynamic Networks
\thanks{A preliminary version of this work will appear in the {\em IEEE
Conference on Decision and Control, 2013.}%
}}

\author{Saghar Hosseini, Airlie Chapman, and Mehran Mesbahi%
\thanks{The research of the authors was supported by the ONR grant N00014-12-1-1002
and AFOSR grant FA9550-12-1-0203-DEF. The authors are with the Department
of Aeronautics and Astronautics, University of Washington, WA 98105.
Emails: \{saghar, airliec, mesbahi\}@uw.edu.%
}}
\maketitle
\begin{abstract}
This paper presents a distributed optimization scheme over a network
of agents in the presence of cost uncertainties and over switching
communication topologies. Inspired by recent advances in distributed
convex optimization, we propose a distributed algorithm based on a
dual sub-gradient averaging. The objective of this algorithm is to
minimize a cost function cooperatively. Furthermore, the algorithm
changes the weights on the communication links in the network to adapt
to varying reliability of neighboring agents. A convergence rate analysis
as a function of the underlying network topology is then presented,
followed by simulation results for representative classes of sensor
networks. \end{abstract}


\section{Introduction}


The past decade has witnessed successful applications of networked
systems in areas ranging from environmental monitoring, robotics,
target recognition, air traffic control, to industrial and manufacturing
automation. By increasing the size and complexity of networked systems,
decentralized optimization schemes are desired for reducing data transmission
rates and ensuring robustness in the presence of local failures. These
methods are particularly relevant when there is a lack of access to
centralized information by individual agents. 

In recent years, there has grown an extensive literature on distributed
convex optimization \cite{Boyd2010,Jakovetic2011,Moskaoyama2007}
and the adaptation and monitoring of the underlying network has become
of increasing interest \cite{Kim2006,Chi2008,Dai2011}. Moreover,
agent\textquoteright s communication range or disturbances may cause
the underlying network topology to change dynamically. In this direction,
a class of distributed sub-gradient algorithms for convex optimization
has been developed \cite{Nedic2009,Lobel2011,Lobel2011a,Duchi2012,Lee2012}.
In these works, local convex cost functions are assumed to be known
while the topology of network is allowed to vary. 

In addition to uncertainties in the network's structure, the environment
can also affect the corresponding cost functions. For such scenarios,
traditional optimization approaches become unsuitable. One approach
to improve the robustness of algorithms for convex optimization is
via stochastic methods \cite{Ram2009,SundharRam2010,Agarwal2011},
where the probability distribution of uncertain variable is known
\textit{a priori}. One such approach has been pursued by Duchi \textit{et
al.} \cite{Duchi2012} who approached this problem using a stochastic
sub-gradient method where the distribution of sub-gradients is known
\textit{a priori}.

Despite its many successes, stochastic optimization-based methods
do not explicitly address the dynamic aspect of the problem in an
unknown environment. Online learning is an extension of stochastic
optimization where the uncertainty in the system is demonstrated by
an \textit{arbitrarily} varying cost function. In particular, at the
time the relevant decision is made the cost function is assumed to
be unknown, without probabilistic assumptions, to the decision-maker.
Such learning algorithms have had a significant impact on modern machine
learning \cite{Zinkevich2003a,Hazan2007a,Xiao2010}. One standard
metric to measure the performance of these online algorithms is called
\textit{regret}. Regret measures the difference between the incurred
cost and the cost of the {\em best fixed decision in hindsight}.
An online algorithm is then declared ``good'' when its regret is
sub-linear.

Distributed online optimization and its applications in multi-agent
systems has not been studied at large by the systems and control community.
Yan \textit{et al.} in \cite{Yan2010} introduced a decentralized
online optimization based on a sub-gradient method in which the agents
interact over a weighted strongly connected directed graph. Considering
an undirected path graph with a fixed-radius neighborhood information
structure, Raginsky \textit{et al.} \cite{Raginsky2011} proposed
an online algorithm for distributed optimization based on sequential
updates, proving a regret bound of $O(\sqrt{T})$. In \cite{Hosseini2013},
we proposed an extension to the work of Duchi \textit{et al.} \cite{Duchi2012}
on distributed optimization with convergence rate of $O(\sqrt{T}\log T)$
to an online setting. In addition, an improved regret bound of $O(\sqrt{T})$
has been derived for strongly connected networks, also highlighting
the dependence of the regret on the connectivity of the underlying
network. 


We note that the aforementioned works do not exploit a dynamic weight
selection procedure to improve the performance of the corresponding
distributed algorithms. In systems and control literature, certain
metrics have been used for designing \textit{adaptive} mechanisms
for networks based on centralized \cite{Ghosh2006,Zelazo2009b,Wan2007}
and distributed \cite{Kim2006,Chapman2010c} strategies. Chapman \textit{et
al.} \cite{Chapman2013cdc} proposed an online distributed algorithm
for re-weighting the network edges in order to dampen the effect of
external disturbances on the system. Dynamic weight selection is also
favorable in the area of sensor networks and distributed estimation
due to power and data rate constraints as well as failure modes of
the inter-sensor communication links \cite{Aldosari2006,Kar2008,Laszka2013}.


In this paper, we consider two types of uncertainties in the networked
systems corresponding to disturbances in cost functions and the network
structure. An adaptive algorithm for distributed optimization over
fixed networks is proposed and further extended to switching graphs.
The main assumption used for implementing this algorithm is that the
local cost function and its sub-gradient are observable at each node
and can be shared with the neighboring nodes in the network. 

The contribution of this paper is threefold. First, we present the
Distributed Weighted Dual Averaging (DWDA) algorithm \cite{Hosseini2013}
for distributed optimization over networks. A distributed dynamic
weight selection method based on an online weighted majority approach
\cite{Littlestone1994,Freund1997} is then embedded in the DWDA algorithm
allowing the weights on the network's edges to adaptively change in
order to optimize the information diffusion in the network. The proposed
algorithm is inspired by the Distributed Dual Averaging (DDA) algorithm
\cite{Duchi2012}. Second, DWDA is applied on switching networks capturing
the uncertainties in communication links. Third, the DWDA is further
extended to Online Distributed Weighted Dual Averaging (Online-DWDA)
algorithm which takes into account the uncertainties in cost functions
and unavailability of reliable statistics on the noise characteristics.
We then proceed to derive regret bounds that highlight the link between
the adaptive weight selection and the Online-DWDA algorithm and can
thus be used to design networks with good regret performance.


The organization of the paper is as follows. The notation and background
on graphs and regret are reviewed in \S\ref{sec:Background-and-Model}.
In \S\ref{sec:Problem-Statement}, the formulation for the distributed
convex optimization problem over networks is presented. This is then
followed by the description for the DWDA algorithm and dynamic weight
selection procedure in \S\ref{sec:Main-Result} and \S\ref{sec:Distributed-Adaptive-Weight};
the convergence analysis of the proposed algorithm is discussed over
switching topologies in \S\ref{sec:Convergence-Analysis}. In \S\ref{sec:Online-Distributed-Optimization},
the distributed convex optimization problem is extended to the online
setting, with applications to networked systems operating in an uncertain
environment. The performance of Online-DWDA is subsequently studied
using the regret analysis. In \S\ref{sec:Simulation results}, we
examine online distributed estimation over sensor networks, demonstrating
the viability of the online approach in distributed estimation. Finally,
\S\ref{sec:Conclusion} provides our concluding remarks and future
directions for utilizing the online framework for system and control
problems.

\section{\label{sec:Background-and-Model}Background and Preliminaries}

We provide a brief background on constructs that will be used in this
paper. For the column vector $v\in\mathbb{R}^{p}$, $v_{i}$ or $\left[v\right]_{i}$
denotes the $i$th element and $e_{i}$ denotes the column vector
which contains all zero entries except $\left[e_{i}\right]_{i}=1$.
The vector of all ones is denoted as $\boldsymbol{1}$. For matrix
$M\in\mathbb{R}^{p\times q}$, $\left[M\right]_{ij}$, or simply $M_{ij}$,
denotes the element in its $i$th row and $j$th column. The family
of probability vectors is denoted by $\Omega$ and contains all non-negative
vectors $\sigma\in\left[0,1\right]^{n}$ such that $\sum\sigma_{i}=1$.
A row stochastic matrix $P$ is a non-negative matrix with rows in
$\Omega$. Moreover, the ergodic coefficient for a stochastic matrix
$Q\in\mathbb{R}^{n\times n}$ is given by 
\begin{equation}
\tau(Q)=1-\min_{i,j\in[n]}\sum_{k=1}^{n}\min\{Q_{ik},Q_{jk}\}.\label{eq:ergodic coef def}
\end{equation}
A time varying matrix is denoted by $P^{t}$ and a (backward) sequence
of time varying matrices is presented by $P^{(t,0)}=P^{t}P^{t-1}\cdots P^{0}.$
For any positive integer $n$, the set $\{1,2,...,n\}$ is denoted
by $\left[n\right]$. The inner product of two vectors $\theta$ and
$w$ is represented by $\left\langle \theta,\phi\right\rangle $.
The 2-norm is signified by $||.||_{2}$, a general norm of vector
is denoted as $\left\Vert \theta\right\Vert $, and its associated
dual norm is defined as $||\theta||_{*}=\underset{||\phi||=1}{\sup}\left\langle \theta,\phi\right\rangle $.
A function $f:\Theta\rightarrow\mathbb{R}$, where $\Theta\subseteq\mathbb{R}^{m}$
for some positive integer $m$, is called $L$-Lipschitz continuous
with respect to the norm $\left\Vert \cdot\right\Vert $ if there
exists a positive constant $L$ for which 
\begin{equation}
|f(\theta)-f(\phi)|\leq L\Vert\theta-\phi\Vert\;\;\mbox{for all}\;\theta,\phi\in\Theta.\label{eq:L-Lipschitz}
\end{equation}

\subsection{\label{sub:Graphs}Graphs}

A succinct way to represent the interactions of dynamic agents, e.g.,
sensors, over a network is through a graph. A weighted directed graph
$\mathcal{G}=\left(V,E,W\right)$ is defined by a node set $V$ with
cardinality $n$, the number of nodes in the graph representing the
agents in the network, and an edge set $E$ comprising of pairs of
nodes which represent the agents' interactions, i.e., agent $i$ affects
agent $j$'s dynamics if there is an edge from $i$ to $j$, denoted
as $\left(i,j\right)\in E$. In addition, a function $W:E\rightarrow\mathbb{R}$
is given that associates a weight $\mbox{\ensuremath{w_{ji}\in W}}$
to every edge $\left(i,j\right)\in E$. Moreover, $\mbox{dist}(i,j)$
denotes the minimum number of edges of any directed path from node
$i$ to node $j$. The time-varying graph topology is presented by
$\mathcal{G}^{t}=\left(V,E^{t},W^{t}\right)$ at time step $t$.
The neighborhood of node $i$ is defined as the set $N_{i}=\{j|(j,i)\in E\}$
and the time-varying neighborhood set is presented by $N_{i}^{t}=\{j|(j,i)\in E^{t}\}$.
The adjacency matrix $A(\mathcal{G})$ is a matrix representation
of $\mathcal{G}$ with $\left[A(\mathcal{G})\right]_{ji}=w_{ji}$
for $\left(i,j\right)\in E$ and $\left[A(\mathcal{G})\right]_{ji}=0$,
otherwise. A graph $\mathcal{G}$ is strongly connected if there exists
a directed path between every pair of distinct vertices. For a graph
$\mathcal{G}$, $d_{i}$ is the weighted in-degree of $i$, defined
as $d_{i}=\sum_{\{j|(j,i)\in E\}}w_{ij}$. In addition, $L(\mathcal{G})=\Delta(\mathcal{G})-A(\mathcal{G})$
is called the graph Laplacian where $\Delta(\mathcal{G})$ is the
diagonal matrix of $d_{i}$'s. Based on the construction of weighted
directed graph Laplacian, every graph $\mathcal{G}$ has a right eigenvector
of $\boldsymbol{1}$ associated with eigenvalue $\lambda=0$~\cite{Mesbahi2010}. 

There are many families of graphs that are often used to model networks
of practical interest. In this paper, we use path and random graphs
for some of our simulations. Specifically, Erd\H{o}s-Rényi random
graphs with edge probability $p$ are constructed by having an edge
$\left(i,j\right)\in E$ in the graph with probability $p$ for all
possible edges. A random tree is a particular realization of a random
graph that is minimally connected and a random $k$-regular graph
is a random graph in which $d_{i}=k$ for all vertices $i\in V$.
In the path graph, $\left(i,j\right)\in E$ if and only if $\left|i-j\right|=1$.

\subsection{Regret\label{sub:Regret background}}

Regret is one measure of performance for learning algorithms. In the
online optimization setting, an algorithm is used to generate a sequence
of decisions $\{x(t)\}_{t=1}^{T}$. The number of iterations is denoted
by $T$ which is unknown to the online player. At each iteration $t$,
after committing to $x(t)$, a previously unknown convex cost function
$f_{t}$ is revealed, and a loss $f_{t}(x(t))$ is incurred. The goal
of the online algorithm is to ensure that the time average of the
difference between the total cost and the cost of the best fixed decision
$x^{*}=\mbox{argmin}\sum_{t=1}^{T}f_{t}(x)$ is small. 
The difference between these two costs over $t=1,2,...,T,$ iterations
is called the regret of the online algorithm, i.e., 
\begin{equation}
R_{T}(x^{*},x)=\sum_{t=1}^{T}\left(f_{t}(x(t))-f_{t}(x^{*})\right).\label{eq:regret_general}
\end{equation}
An algorithm performs well if its regret is sub-linear as a function
of $T$, i.e. $\lim_{T\rightarrow\infty}R_{T}/T=0$. This implies
that on average, the algorithm performs as well as the best fixed
strategy in hindsight independent of the adversary's moves and environmental
uncertainties. Further discussion on online algorithms and their regret
analysis can be found in \cite{Shalev-Shwartz2011b,SebastienBubeck2011,Hazan2011}.

The general definition of regret is presented in \eqref{eq:regret_general}
for a single decision-making unit. 
In order to analyze the performance of {\em distributed online algorithms}
two variations of the notion of regret are introduced. First is the
regret due to agent $i$'s decision, 
\begin{equation}
R_{T}(x^{*},x_{i})=\sum_{t=1}^{T}\left(f_{t}(x_{i}(t))-f_{t}(x^{*})\right),\label{eq:regret_distributed}
\end{equation}
which is the cumulative penalty agent $i$ incurs due to its local
decisions $\{x_{i}(t)\}_{t=1}^{T}$ on the global cost sequence $\left\{ f_{t}\right\} $.
Second is the regret based on the {\em running average} of the
decisions $\{x_{i}(t)\}_{t=1}^{T}$, 
\begin{equation}
R_{T}(x^{*},\widetilde{x}_{i})=\sum_{t=1}^{T}\left(f_{t}(\widetilde{x}_{i}(t))-f_{t}(x^{*})\right),\label{eq:regret_general-average}
\end{equation}
where $\widetilde{x}_{i}(T)=\frac{1}{T}\sum_{t=1}^{T}x_{i}(t).$

\section{Problem Statement\label{sec:Problem-Statement}}

In this section a distributed decision process is considered in which
a large number of agents cooperatively optimize a global objective
function over the network denoted by $\mathcal{G}=(V,E,W)$. 

The global objective to be minimized is
\begin{equation}
f(x)=\frac{1}{n}\sum_{i=1}^{n}f_{i}(x)\quad\mbox{subject to}\;\; x\in\chi,\label{eq:objective fun-1}
\end{equation}
where $f_{i}(x):\mathbb{\mathbb{R}}^{d}\rightarrow\mathbb{R}$ is
a convex cost function associated with agent $i\in V$ and $\chi\subseteq\mathbb{\mathbb{R}}^{d}$
is a closed convex set. The global optimization problem will be solved
\textit{locally} by each agent $i$ via the local decision variable
$x_{i}\in\chi$.

\section{Weighted Dual Averaging\label{sec:Main-Result}}

In order to solve the optimization problem \eqref{eq:objective fun-1},
we adapt Nesterov's dual averaging algorithms \cite{Nesterov2007}
and our preliminary results on the Distributed Weighted Dual Averaging
(DWDA) algorithm \cite{Hosseini2013}, that in turn is inspired by
\cite{Duchi2012}. The DWDA algorithm sequentially updates the \textit{local}
$x_{i}(t)$ and a \textit{working} variable $y_{i}(t)$ for each agent
$i$. The update itself is based on a provided local sub-gradient
of the loss $f_{i}(x_{i}(t))$ denoted as $g_{i}(t)$. The centralized
form of the dual averaging algorithm appears as a sub-gradient decent
method followed by a projection step onto the constraint set $\chi$,
specifically, 
\begin{equation}
y(t+1)=y(t)+g(t),\label{eq:gradient decent}
\end{equation}
where $g(t)\in\partial f(x(t))$. Then 
\begin{equation}
x(t+1)=\Pi_{\chi}^{\psi}\left(y(t+1),\alpha(t)\right),\label{eq:projection step}
\end{equation}
where $\Pi_{\chi}^{\psi}(\cdot)$ is a regularized projection onto
$\chi$, to be formally defined shortly. The Distributed Weighted
Dual Averaging (DWDA) algorithm is presented in Algorithm \ref{alg:DWDA}.
The projection function used in this algorithm is defined as 
\begin{equation}
\Pi_{\chi}^{\psi}(y(t),\alpha(t))=\arg\min_{x\in\chi}\left\{ \langle y(t),x\rangle+\frac{1}{\alpha(t)}\psi(x)\right\} ,\label{eq:projection fun}
\end{equation}
where $\alpha(t)$ is a non-increasing sequence of positive functions
and $\psi(x):\chi\rightarrow\mathbb{R}$ is a proximal function. The
standard dual averaging algorithm uses proximal function $\psi(x)$
to avoid undesirable oscillations in the projection step. Without
loss of generality, $\psi$ is assumed to be strongly convex with
respect to $\Vert.\Vert,$ $\psi(x)\geq0$, and $\psi(0)=0$.

\begin{algorithm}
\small
\SetAlgoLined 
\For {$t=1$ \KwTo $T$}{
	Evaluate $f(t)=\{f_{i}(t);\;\mbox{for all}\; i=1,...,n\}$\\
	\ForEach {Agent $i$}{
  Compute subgradient $g_i(t)\in \partial f_{i}({x}_{i}(t))$\\
     $y_i(t+1) = \sum_{j\in N(i)} P_{ji}(t)y_j(t) + g_i(t)$\ \label{eq:z_i update}

     ${x}_i(t+1) = \prod_{\chi}^{\psi}(y_i(t+1),\alpha(t))$\ \label{eq:x_i update}

     $\widetilde{x}_i(t+1)=\frac{1}{t+1}\sum_{s=1}^{t+1} {x}_i(s)$\ \label{eq:x_avg update}
 }}
\caption{Distributed Weighted Dual Averaging (DWDA)} 
\label{alg:DWDA}
\end{algorithm}

The distributed algorithm can be considered as an approximated sub-gradient
descent. The approximation is attained by an agent via a convex combination
of local sub-gradients provided by its neighbors. This operation can
be represented compactly as a stochastic matrix $P\in\mathbb{R}^{n\times n}$
which preserves the zero structure of the Laplacian matrix $L(\mathcal{G})$.
It is clear that for all agents to have access to each cost function
$f_{i}$ there must be a path from every agent $i$ to every other
agent. Consequently, a minimum requirement on the underlying network
is that it must be strongly connected. The distributed dynamic weight
selection procedure presented in the following section constructs
a row stochastic matrix $P$ of the required form that is associated
with a weighted directed graph.

\section{Distributed Dynamic Weight Selection\label{sec:Distributed-Adaptive-Weight}}

In this section, we propose an adaptation scheme for the network weight
selection in order to improve the information diffusion in line \ref{eq:z_i update}
of Algorithm \ref{alg:DWDA} such that the communication matrix $P$
is a row stochastic matrix with positive diagonal elements. In this
proposed distributed algorithm, each agent $i\in[n]$ estimates its
loss function via a convex combination of loss functions available
to it by its neighboring agents. This convex combination is specified
by weights $w_{ij}$'s on each edge $\left(j,i\right)\in E$ and $w_{ii}$,
respectively. The edge re-weighting problem parallels the Weighted
Majority (WM) algorithm \cite{Littlestone1994}. The context of the
WM algorithm is the presence of $|N_{i}|+1$ experts, and the associated
cost $h_{j}(t)$ assigning a loss value to expert $j\in\{N_{i},i\}$,
where $0\leq h_{j}(t)\leq1$. At each time-step, agent $i$ selects
a probability distribution $q(t)$ over the $|N_{i}|+1$ experts,
i.e., $q(t)\in\mathcal{Q}_{i}=\left\{ q\in\Omega|q_{j}=0\mbox{ for }j\notin\{N_{i},i\}\right\} $,
in order to minimize $l_{t,i}(q)=\sum_{j\in\{N(i),i\}}q_{j}(t)h_{j}(t).$
The regret for agent $i$ for the WM algorithm is then defined as
\begin{equation}
L_{T}(q^{*},q)=\sum_{t=1}^{T}l_{t,i}(q)-\sum_{t=1}^{T}l_{t,i}(q^{*}),\label{eq:WM regret}
\end{equation}
where 
\[
q^{*}=\mbox{argmin}_{q\in\mathcal{Q}_{i}^{'}}\sum_{t=1}^{T}f_{t}(q_{j})
\]
 with $\mathcal{Q}_{i}^{'}=\left\{ q\in\mathcal{Q}_{i}|q_{j}\in\left\{ 0,1\right\} \right\} $,
and the best fixed strategy $q^{*}$ is the best expert $j\in\{N_{i},i\}$
in hindsight. Consequently, \eqref{eq:WM regret} is of the same form
as \eqref{eq:regret_general}.

The general form of the WM algorithm is presented in \cite{Freund1997}
as the Online Allocation (OA) algorithm that is applicable to any
bounded loss function over general decision and outcome spaces. Based
on the OA algorithm, the regret for each agent $i\in[n]$ is bounded
as 
\begin{equation}
L_{T}\leq M\left(\sqrt{2T\ln\left(|N_{i}|+1\right)}+\ln\left(|N_{i}|+1\right)\right),\label{eq:OA regret}
\end{equation}
where $M$ is the upper bound on the loss function $h_{j}(t)$ for
all $j\in\{N_{i},i\}$. Since the regret in \eqref{eq:OA regret}
is sub-linear over time, the weight allocation performs as well as
the best strategy in hindsight. 

A Distributed Online Allocation (DOA) algorithm is proposed based
on the OA algorithm where agent $i\in[n]$ specifies the weights $w_{ij}$'s
associated to each edge $\left(j,i\right)\in E$ as well as the weight
on the self-loop $w_{ii}$. The DOA algorithm presented in Algorithm
\ref{alg:DOA} is embedded in Algorithm \ref{alg:DWDA} at each iteration. 

\begin{algorithm}
\small
\SetAlgoLined 
\ForEach {Agent $i$}{
Choose $\beta \in [0,1]$\ and initial weight vector $\mathbf{w}_i(t)=\mathbf{1} $\\
Let $q(t)=\frac{1}{n} \mathbf{w}_i(t)$\\
\For {$t=1$ \KwTo $T$}{
Adversary reveals $f(t)=\{f_{j}(t);\;\mbox{for all}\; j\in \{N(i),i\}\}$\\
Suffer loss $l_{t,i} = \sum_{j\in \{N(i),i\}} q_{j}(t)f_{j}(t)$\\
\ForEach {Agent $j\in [n]$}{
	\eIf {$j\in \{N(i),i\}$}{
     $w_{ij}(t+1) = w_{ij}(t) \beta ^{f_{j}(t)}$ \label{eq:w_ij update}
    }{
     $w_{ij}(t+1) = w_{ij}(t)$\
 }}
$q(t+1) = \frac{\mathbf{w}_i(t+1)}{\sum_{j\in \{N(i),i\}} w_{ij}(t+1)} $\ \label{eq:q_i update}\\
}}
\caption{Distributed Online Allocation (DOA)}  
\label{alg:DOA}
\end{algorithm}

In the distributed optimization process considered, each agent decides
on the weights associated with the information \textit{received} from
its neighboring agents. This information is based on the neighbor's
local loss function. Intuitively, the algorithm places more weight
on the link associated with the neighboring agent that has a higher
confidence in its decision. The positive diagonal entries represent
the self-confidence of each agent and is updated based on the local
loss. 

In networks with fixed topology, the communication matrix $P(t)$
preserves its zero structure for all time $t$. In addition, the non-zero
elements in each row of the communication matrix $P$ is specified
by line \ref{eq:q_i update} in Algorithm \ref{alg:DOA}: 
\begin{equation}
P_{ij}(t)=\begin{cases}
q_{j}(t) & \mbox{for}\; j\in\{N_{i},i\}\\
0 & \textrm{otherwise}
\end{cases}.\label{eq:P row stochastic}
\end{equation}
Since for each agent $i$, $q(t)\in\mathcal{Q}_{i}$ is a probability
distribution, the communication matrix $P(t)$ will be row stochastic
at every time step. The weighted graph Laplacian can then be formed
as 
\begin{equation}
L(\mathcal{G}(t))=I-P(t).\label{eq:weighted graph laplacian}
\end{equation}
In addition, note that since the graph is strongly connected, the
communication matrix $P$ is 1-irreducible (\cite{Wu2005}; Corollary
4) and given positive diagonal elements, it is in-decomposable and
aperiodic (SIA). In the following section, the DOA algorithm is extended
to construct a row stochastic communication matrix $P(t)$ for directed
switching graphs with time-varying edge sets.

\subsection{Switching Topologies\label{sub:Switching-Topologies}}

The network topology may change dynamically due to disturbances or
communication range limitations. In this section we apply the dynamic
weight selection procedure discussed in \S \ref{sec:Distributed-Adaptive-Weight}
to switching topologies. In this paper we assume that the union of
directed topologies $\mathcal{G}^{\cup_{i=0}^{\delta-1}}=\bigcup_{i=0}^{\delta-1}\mathcal{G}^{i}$
over some fixed uniform intervals $\delta$, with $\delta\geq1$ a
positive integer, is strongly connected. We note that the communication
matrix of $\mathcal{G}^{\cup_{i=0}^{\delta-1}}$ can be presented
as 
\begin{equation}
P^{\cup_{i=0}^{\delta-1}}=P^{0}+P^{1}+\dots+P^{\delta-1}.\label{eq:union communication matrix}
\end{equation}
Thus, each row of the communication matrix $P$ over switching topologies
is specified by line \ref{eq:q_i update} in Algorithm \ref{alg:DOA}
as

\begin{equation}
P_{ij}(t)=\begin{cases}
q_{j}(t)/\left(\sum_{k\in\{N_{i}^{t},i\}}q_{k}(t)\right) & \mbox{for}\; j\in\{N_{i}^{t},i\}\\
0 & \textrm{otherwise}
\end{cases},\label{eq:P row stochastic switching}
\end{equation}
where $q(t)\in Q_{i}^{t}=\left\{ q\in\Omega|q_{j}=0\mbox{ for }j\notin\{N_{i}^{t},i\}\right\} $
and is a probability distribution. Note that the communication matrix
$P(t)$ will be row stochastic at every time step and thus the weighted
graph Laplacian is the same as in \eqref{eq:weighted graph laplacian}.
Since the graph $\mathcal{G}^{\cup_{i=0}^{\delta-1}}$ is strongly
connected and $P^{\cup_{i=0}^{\delta-1}}$ has positive diagonal elements,
the communication matrix $P^{\cup_{i=0}^{\delta-1}}$ is SIA (\cite{Wu2005};
Corollary 4). These properties of communication matrices will be subsequently
employed in the convergence analysis of the DWDA algorithm.

\section{Convergence Analysis\label{sec:Convergence-Analysis}}

Before presenting the convergence analysis of the distributed optimization
algorithm, a few preliminary remarks and assumptions are in order.
We assume that each convex function $f_{i}$ is positive and $L$-Lipschitz
with respect to $\Vert.\Vert$. Assuming that $\left\{ P^{t}\right\} $
is SIA, there exists a vector $\pi\in\Omega$~\cite{Anthonisse1977},
such that 
\begin{equation}
\pi_{j}=\sum_{i=1}^{n}\pi_{i}P_{ij}^{t}\;\;\mbox{for all}\; t\in[T],\label{eq: weight factors}
\end{equation}
where $\pi_{i}$ is referred to as the weighting factor for agent
$i$.

In order to take advantage of the properties of the standard weighted
dual averaging in our regret analysis, the sequences $\bar{y}(t)$
and $\bar{g}(t)$ are defined as 
\begin{align}
\bar{y}(t) & =\sum_{i=1}^{n}\pi_{i}y_{i}(t),\mbox{ and }\bar{g}(t)=\sum_{i=1}^{n}\pi_{i}g_{i}(t),\label{eq:average z,g}
\end{align}
signifying the (network-level) weighted average of dual variables
and subgradients in the DWDA algorithm, respectively. Therefore, based
on \eqref{eq: weight factors} and \eqref{eq:average z,g}, 
\begin{align}
\bar{y}(t+1) & =\sum_{i=1}^{n}\pi_{i}\left\{ \sum_{j=1}^{n}P_{ij}^{t}y_{j}(t)+g_{i}(t)\right\} \nonumber \\
 & =\sum_{\begin{array}{c}
j=1\end{array}}^{n}y_{j}(t)\sum_{i=1}^{n}\pi_{i}P_{ij}^{t}+\bar{g}(t)\nonumber \\
 & =\sum_{\begin{array}{c}
j=1\end{array}}^{n}y_{j}(t)\pi_{j}+\bar{g}(t)\nonumber \\
 & =\bar{y}(t)+\bar{g}(t),\label{eq:z_bar update}
\end{align}
which is analogous to the dual averaging update \eqref{eq:gradient decent}.
Thus, the following update rule is introduced which is analogous to
the standard dual averaging algorithm projection step \eqref{eq:projection step},
where the primal variable is updated as 
\begin{equation}
\phi(t+1)=\Pi_{\chi}^{\psi}(\bar{y}(t+1),\alpha(t)).\label{eq:y update}
\end{equation}
The performance analysis of the distributed optimization and adaptive
weight selection can now be presented. 

The following result by Duchi \textit{et al.} implies that after $T$
iterations of Algorithm \ref{alg:DWDA}, each agent's error in the
evaluation of total cost is bounded by the error due to Dual-Averaging
method.
\begin{thm}
\cite{Duchi2012}\label{thm:offline Regret(xi)} Given the sequences
$x_{i}(t)$ and $y_{i}(t)$ generated by lines \textup{\ref{eq:z_i update}}
and \textup{\ref{eq:x_i update}} in \textup{Algorithm \ref{alg:DWDA}},
for all $i\in[n]$ with proximal function $\psi$ and $\alpha(t)>0$,
we have 
\begin{align}
 & \frac{1}{T}\sum_{t=1}^{T}f(x_{i}(t))-f(x^{*})\leq\frac{L^{2}}{2}\sum_{t=1}^{T}\alpha(t-1)+\frac{1}{\alpha(T)}\psi(x^{*})\nonumber \\
 & +\frac{L}{T}\sum_{t=1}^{T}\alpha(t)(\Vert\bar{y}(t)-y_{i}(t)\Vert_{*}+\frac{2}{n}\sum_{i=1}^{n}\Vert\bar{y}(t)-y_{i}(t)\Vert_{*}).\label{eq:offline f-f*}
\end{align}

\end{thm}
The last two terms on the right hand side of \eqref{eq:offline f-f*}
represent the error due to the network which is defined as the deviation
of local dual variable $y_{i}$ from the weighted average of dual
variables $\bar{y}$ over the network. Lemma \ref{lem:z-z_i-1} in
the Appendix imposes an upper bound on the effect of network topology
associated with $\Vert\bar{y}(t)-y_{i}(t)\Vert_{*}$ as 
\begin{equation}
\Vert\bar{y}(t)-y_{i}(t)\Vert_{*}\leq L\sum_{k=0}^{t-2}\sum_{j=1}^{n}\left|P_{ij}^{(t-1,k+1)}-\pi_{j}\right|+2L.\label{eq:netwrok effect}
\end{equation}

Inequality \eqref{eq:netwrok effect} highlights the importance of
the underlying network topology through the communication matrix $P^{t}$
and its products. Note that the network effect is analogues to the
consensus problems \cite{book}. Therefore we proceed to extend the
distributed optimization algorithm to switching graphs in the following
subsection and provide a sub-linear convergence rate for the DWDA
algorithm.

\subsection{Switching Topologies}

In this section we employ the weak ergodicity of inhomogeneous Markov
chains to reason about the convergence of the DWDA algorithm. This
property implies that the product of stochastic SIA matrices converges
exponentially to a rank-one matrix of the form $\mathbf{1}\pi^{T}$
as $t\rightarrow\infty$, where $\pi\in\Omega$. Applying the following
result from \cite{Jadbabaie2003}, it thus suffices to show that $P^{(\delta-1,0)}$
is SIA.
\begin{lem}
\cite{Jadbabaie2003}\label{lem: union of P} Let $m\geq1$ be a positive
integer and $P^{\tau}$be non-negative matrices with positive diagonal
elements for $\tau=0,1,\dots,m$. Then, $P^{m}P^{m-1}\dots P^{0}\geq\mu^{m}(P^{0}+P^{1}+\dots+P^{m}),$
where $\mu>0$ is specified by the diagonal elements of matrices $P^{\tau}$
for all $\tau=0,1,\dots,m$. 
\end{lem}
From Lemma \ref{lem: union of P} and \eqref{eq:union communication matrix},
$P^{(\delta-1,0)}$ is bounded below by an SIA matrix. Moreover, we
note that $P^{(\delta-1,0)}$ is also a stochastic matrix, thus it
must be SIA. Therefore, based on the weak ergodicity of inhomogeneous
Markov chains, the product 
\[
P^{(k\delta-1,(k-1)\delta)}\dots P^{(2\delta-1,\delta)}P^{(\delta-1,0)}
\]
converges exponentially to a rank-one matrix of the form $\mathbf{1}\pi^{T}$
as $t\rightarrow\infty$, and based on Theorem 1 of \cite{Anthonisse1977},
we have

\begin{equation}
\left|P^{(k\delta-1,0)}-\pi_{j}\right|\leq\gamma^{\left\lfloor \frac{k}{\nu}\right\rfloor },\label{eq:perron-Frobenious-2}
\end{equation}
where 
\begin{equation}
\gamma=\max_{\nu\geq1}\left\{ \tau(P^{(\delta\nu-1,0)})<1\right\} .\label{eq:gamma-2}
\end{equation}
Note that the maximization is over all realizations of the sequence
$P^{(\delta\nu-1,0)}$ and $\nu$ is bounded as stated in the following
proposition.
\begin{prop}
\label{prop: switching nu}Consider a set $\mathcal{P}$ of stochastic
matrices with positive diagonal elements, representing arbitrarily
strongly connected topologies over $n$ nodes, i.e., $P^{k}\in\mathcal{P}$
for all positive integers $k$. Then, there exists an integer $\nu$,
$1\leq\nu\leq n-1$, for which if the sequence $Q=P^{(m+\nu-2,m)}$
of matrices in $\mathcal{P}$ is not scrambling, $P^{m+\nu-1}Q$ is
scrambling. \end{prop}
\begin{IEEEproof}
Let $Q_{1}=P^{m}$ and $Q_{2}=P^{m+1}P^{m}.$ Then every entry of
$Q_{2}$ is represented as 
\[
\left[Q_{2}\right]_{ij}=\sum_{k=1}^{n}\left[P^{m+1}\right]_{ik}\left[P^{m}\right]_{kj}.
\]
Let $E^{\cup_{i=m}^{m+1}}$ represent the edge set of the union of
directed graphs associated with $P^{m}$ and $P^{m+1}$. Since $\left[P^{m}\right]_{ii}>0$
for all $i\in[n]$ and integer $m\geq1$, the entry $\left[Q_{2}\right]_{ij}$
is positive if $\left(j,i\right)\in E^{\cup_{i=m}^{m+1}}$ or if there
exists a node $k\in[n]$ such that $(j,k)\in E^{m}$ and $(k,i)\in E^{m+1}$.
Thus, the corresponding zero entry of $Q_{1}$ that has one of the
aforementioned properties will be positive in $Q_{2}.$ By induction,
it follows that the entry of $\left[Q_{\nu}\right]_{ij}$ will be
positive if $\left(j,i\right)\in E^{\cup_{i=m}^{m+\nu-1}}$, or if
there exists a set of nodes $\{k_{m},k_{m+1},\dots,k_{m+\nu-2}\}$
such that $\{(j,k_{m})$, $(k_{m},k_{m+1})$, ..., $(k_{m+\nu-2},i)\}\in E^{m}\times E^{m+1}\times\dots\times E^{m+\nu-1}$.
Therefore, for each row $i$ of $Q_{\nu}$, all entries will be positive
when 
\begin{equation}
\nu_{i}=\max_{j\in\left[n\right],m\geq1}\left\{ \mbox{dist}(j,i)\,\mbox{for }\mathcal{G}^{m}=\left(V,E^{m},W\right)\right\} .\label{eq:v_i_switching}
\end{equation}
Note that the maximization in \eqref{eq:v_i_switching} is over all
possible strongly connected graphs with the directed cycle graph representing
the worst case with $\nu_{i}=n-1$. Since every element of any row
of the sequence $Q=P^{(m+\nu-2,m)}$ of matrices in $\mathcal{P}$
is positive, the matrix $Q$ is scrambling and thus $1\leq\nu\leq n-1$.
\end{IEEEproof}
We also note that the fixed topology is a special case of switching
graphs with $\delta=1$ in \eqref{eq:union communication matrix}.
Moreover, Proposition \ref{prop: shortest path} of the Appendix presents
a less conservative bound on $\nu$ for fixed topologies. 

Now, we can state the following theorem for the rate of convergence
of DWDA over switching graphs.
\begin{thm}
\label{thm:offline_conv_rate}Given the sequences $x_{i}(t)$ and
$y_{i}(t)$ generated by lines \textup{\ref{eq:z_i update}} and \textup{\ref{eq:x_i update}}
in \textup{Algorithm \ref{alg:DWDA}}, for all $i\in[n]$ with $\psi(x^{*})\leq R^{2}$
and $\alpha(t)=k/\sqrt{t}$, we have 
\begin{align}
 & \frac{1}{T}\sum_{t=1}^{T}f(x_{i}(t))-f(x^{*})\leq\nonumber \\
 & \left(\frac{R^{2}}{k}+kL^{2}\left(\frac{6n}{1-\gamma}+6n\delta\nu+1\right)\right)\frac{1}{\sqrt{T}},\label{eq:convergence rate}
\end{align}
where $\gamma<1$ is a function of the ergodicity of the communication
matrix (see \eqref{eq:gamma-2}) while $\nu$ is a measure of network
connectivity and is bounded by the diameter of the graph $\mathcal{G}^{t}$
(see also Proposition \ref{prop: switching nu}). In addition, $k>0$
is an arbitrary constant and $\delta\geq1$ is a positive integer
as presented in \eqref{eq:union communication matrix}.\end{thm}
\begin{IEEEproof}
Based on \eqref{eq:netwrok effect} and \eqref{eq:perron-Frobenious-2}
we have 
\begin{equation}
\Vert\bar{y}(t)-y_{i}(t)\Vert_{*}\leq nL\sum_{k=1}^{t-1}\gamma^{k}+nL(\delta\nu-1)+2L,\label{eq:network effect 2}
\end{equation}
and since $\gamma<1$, \eqref{eq:network effect 2} is further bounded
as 
\begin{equation}
\Vert\bar{y}(t)-y_{i}(t)\Vert_{*}\leq nL\left(\frac{1}{1-\gamma}+\delta\nu-1\right)+2L.\label{eq:network effedct 3}
\end{equation}
Therefore, the integral test on $\alpha(t)=k/\sqrt{t}$ provides a
bound%
\footnote{Note that $\sum_{t=1}^{T}\frac{k}{\sqrt{t}}\leq2k\sqrt{T}-k$.%
} on the first and last terms in \eqref{eq:convergence rate} as
\begin{align*}
\frac{1}{T}\sum_{t=1}^{T}f(x_{i}(t))-f(x^{*}) & \leq\frac{kL^{2}}{\sqrt{T}}+\frac{\psi(\theta^{\star})}{k\sqrt{T}}+\\
 & \;\;\;\;\frac{6kL^{2}}{\sqrt{T}}\left(\frac{n}{1-\gamma}+n\delta\nu+1\right).
\end{align*}
Given $\psi(x^{*})\leq R^{2},$ the statement of the theorem now follows. 
\end{IEEEproof}
Theorem \ref{thm:offline_conv_rate} states that Algorithm \ref{alg:DWDA}
performs ``well'' as it exhibits a sub-linear convergence rate.
It also highlights the importance of the underlying network topology
through the parameters $\gamma$ and $\nu$. In particular, $\nu$
corresponds to the diameter of the graph as expressed in Proposition
\ref{prop: switching nu} and $\gamma$ is proportional to the ergodic
coefficient $\tau(P^{t})$ of the communication matrix $P^{t}$ as
formed by Algorithm \ref{alg:DWDA}. The ergodic coefficient bounds
the second largest eigenvalue of $P^{t}$, $\lambda_{2}(P^{t})$,
as $|\lambda_{2}(P^{t})|\leq\tau(P^{t})<1$. Thus, based on \eqref{eq:weighted graph laplacian},
$1-\lambda_{2}(P^{t})=\lambda_{n-1}(\mathcal{G}^{t})$ where $\lambda_{n-1}(\mathcal{G}^{t})$
is the second smallest eigenvalue of the weighted graph Laplacian
$L(\mathcal{G}^{t})$ and a well known measure of network connectivity.
Consequently, high network connectivity promotes good performance
of the proposed algorithm. 

In the following section we study the effect of the proposed dynamic
weight selection on the network connectivity and the convergence rate
\eqref{eq:convergence rate}.

\subsection{Adaptive Weight Selection\label{sub:Adaptive-Weight-Selection}}

In this section, we show that embedding Algorithm \ref{alg:DOA} within
Algorithm \ref{alg:DWDA} improves the network information flow and
the speed of convergence in \eqref{eq:convergence rate}. To this
end the following result provides a bound on the ergodic coefficient.
\begin{thm}
\label{thm: ergodic coeff UB switching}Suppose the sequence $q_{j}(t)$
generated by line \textup{\ref{eq:q_i update}} of Algorithm\textup{
\ref{alg:DOA} and }the communication matrices $P^{t}$ are constructed
by \eqref{eq:P row stochastic switching}. Then\textup{ 
\[
\tau(P^{(\delta\nu-1,0)})\leq1-\frac{n}{(\max_{i\in[n],t\in[\delta\nu-1]}|N_{i}^{t}|+1)^{\delta\nu}},
\]
}where $\nu$ is a measure of network connectivity and is bounded
by the diameter of the graph $\mathcal{G}^{t}$ (see also Proposition
\ref{prop: switching nu}). In addition, $\delta\geq1$ is a positive
integer as presented in \eqref{eq:union communication matrix}.\end{thm}
\begin{IEEEproof}
Based on line \ref{eq:q_i update} of Algorithm \ref{alg:DOA}, we
have that for all $k\in\{N_{i},i\}^{t}$, 
\begin{equation}
q_{k}(t)=\frac{\beta^{-\left(\sum_{s=1}^{r_{k,i}(t-1)}f_{k}(s)\right)}}{\sum_{j\in\{N_{i}^{t-1},i\}}\beta^{-\left(\sum_{s=1}^{r_{j,i(t-1)}}f_{j}(s)\right)}},\label{eq:q_k switching}
\end{equation}
where $r_{ji}(t)$ represents the number of communication rounds through
the directed edge $(j,i)$ up to time $t$. Subsequently \eqref{eq:P row stochastic switching}
and \eqref{eq:q_k switching} imply 
\[
P_{ik}^{t}=\begin{cases}
\frac{\beta^{-\left(\sum_{s=1}^{r_{k,i}(t-1)}f_{k}(s)\right)}}{\sum_{j\in\{N_{i}^{t},i\}}\beta^{-\left(\sum_{s=1}^{r_{j,i(t-1)}}f_{j}(s)\right)}} & \left(k,i\right)\in E^{t},\\
0 & \mbox{otherwise}.
\end{cases}
\]
Since $\min_{j\in[n]}\sum_{s=1}^{r_{j,i}(t-1)}f_{j}(s)\leq\sum_{s=1}^{r_{k,i}(t-1)}f_{k}(s)$,
we have 
\[
\beta^{-\left(\sum_{s=1}^{r_{k,i}(t-1)}f_{k}(s)\right)}\leq\beta^{-\left(\min_{j\in[n]}\sum_{s=1}^{r_{j,i}(t-1)}f_{j}(s)\right)}
\]
 and one can bound $P_{ik}^{t}$ from below for all $\left(k,i\right)\in E$
as 
\begin{equation}
P_{ik}^{t}\geq\frac{1}{(|N_{i}^{t}|+1)}\beta^{-C_{k}(t-1)},\label{eq:P_ij(t) LB-1}
\end{equation}
where $C_{k}(t)=\sum_{s=1}^{r_{k,i}(t)}f_{k}(s)-\min_{j\in[n]}\sum_{s=1}^{r_{j,i}(t)}f_{j}(s)$.
Since $C_{k}(t-1)\geq0$ for all $k\in[n]$ and $\beta\in[0,1]$,
we have 
\[
P_{ij}^{t}\geq\frac{1}{(|N_{i}^{t}|+1)},
\]
 for all $t\in[T]$ and subsequently

\[
P_{ij}^{(\delta\nu-1,0)}\geq\frac{1}{(\max_{i\in[n],t\in[\delta\nu-1]}|N_{i}^{t}|+1)^{\delta\nu}}.
\]
Based on \eqref{eq:ergodic coef def}, the statement of the theorem
now follows.
\end{IEEEproof}
Theorem \ref{thm: ergodic coeff UB switching} in conjunction with
\eqref{eq:gamma-2} imply 
\[
\gamma\leq1-\frac{n}{(\max_{i\in[n],t\in[\delta\nu-1]}|N_{i}^{t}|+1)^{\delta\nu}},
\]
which proves to be a conservative bound as the DOA algorithm leads
to a tighter upper bound capturing the performance of agents. In other
words, based on \eqref{eq:P_ij(t) LB-1}, we can show that 
\[
\gamma\leq1-\frac{n\beta^{-J_{\delta\nu}}}{(\max_{i\in[n]}|N_{i}|+1)^{\delta\nu}},
\]
where $J_{\delta(n-1)}=C_{i}(1)+C_{j}(2)+...+C_{k}(\delta\nu)$ and
$i,j,\dots,k\in[n]$. In addition, we know that $\beta^{-J_{\delta\nu}}>1$
and $J_{\delta\nu}$ is an increasing sequence. If the agents are
not performing well, $J_{\delta\nu}$ will increase and subsequently
$\tau(P^{(\delta\nu,0)})$ will decrease which suggests that the DOA
algorithm mitigates the effect of the network topology in \eqref{eq:convergence rate}.
Moreover, Theorem \ref{thm: ergodic coeff UB switching} implies that
the DWDA algorithm performs well for certain types of graphs such
as $k$-regular and expander graphs where the maximum number of neighbors
can be bounded.

\section{Online Distributed Optimization\label{sec:Online-Distributed-Optimization}}

We now consider the effect of uncertainties in the environment on
distributed decision processes where the global objective is to minimize
\begin{equation}
f(x)=\frac{1}{n}\sum_{i=1}^{n}f_{t,i}(x)\quad\mbox{subject to}\;\; x\in\chi,\label{eq:objective fun}
\end{equation}
where $f_{t,i}:\mathbb{\mathbb{R}}^{d}\rightarrow\mathbb{R}$ is a
convex cost function associated with agent $i\in[n]$, assumed to
be revealed to the agent only after the agent commits to the decision
$x(t)$. In other words, the function $f_{t,i}$ is allowed to change
over time in an unpredictable manner due to modeling errors and uncertainties
in the environment. The optimization variable $x_{i}\in\mathbb{R}^{d}$
belongs to a closed convex set $\chi\subseteq\mathbb{\mathbb{R}}^{d}$
and represents the local decision made by agent $i$. Furthermore,
the online-DWDA scheme is analogous to the DWDA presented in Algorithm
\ref{alg:DWDA}. The regret analysis is presented in the following
result quantifying the performance of the proposed algorithm. 
\begin{thm}
\label{thm:Regret(xi)}Given the sequences $x_{i}(t)$ and $y_{i}(t)$
generated by lines \textup{\ref{eq:z_i update}} and \textup{\ref{eq:x_i update}}
in \textup{Algorithm \ref{alg:DWDA}}, for all $i\in[n]$ with $\psi(x^{*})\leq R^{2}$
and $\alpha(t)=k/\sqrt{t}$, we have 
\begin{equation}
R_{T}(x^{*},x_{i})\leq\left(\frac{R^{2}}{k}+kL^{2}\left(\frac{6n}{1-\gamma}+6n\delta\nu+1\right)\right)\sqrt{T},\label{eq:regret bound}
\end{equation}
where $\gamma$ is a function of the ergodicity of the communication
matrix (see \eqref{eq:gamma-2}) while $\nu$ is a measure of network
connectivity and is bounded by the diameter of the graph $\mathcal{G}^{t}$
(see also Proposition \ref{prop: shortest path}). In \eqref{eq:regret bound},
$k>0$ is an arbitrary constant and $\delta\geq1$ is a positive integer
at which the union of directed topologies $\mathcal{G}^{\cup_{i=0}^{\delta-1}}=\bigcup_{i=0}^{\delta-1}\mathcal{G}^{i}$
over some fixed uniform intervals $\delta$ is strongly connected.\end{thm}
\begin{IEEEproof}
Consider an arbitrary fixed decision $x^{*}\in\chi$ and a sequence
$\phi(t)$ generated by \eqref{eq:y update}. From the $L$-Lipschitz
continuity of $f_{t,i}$'s and the definition of regret in \eqref{eq:regret_distributed},
the regret is bounded as 
\begin{equation}
R_{T}(x^{*},x_{i})\leq\sum_{t=1}^{T}\left(f_{t}(\phi(t))-f_{t}(x^{*})+L\Vert x_{i}(t)-\phi(t)\Vert\right).\label{eq: th1_1}
\end{equation}
Note that we can reformulate the first term on the right hand side
of \eqref{eq: th1_1} as 
\begin{align}
f_{t}(\phi(t))-f_{t}(x^{*}) & =\left(\frac{1}{n}\sum_{i=1}^{n}f_{t,i}(x_{i}(t))-f_{t}(x^{*})\right)\nonumber \\
 & +\left(\frac{1}{n}\sum_{i=1}^{n}\left[f_{t,i}(\phi(t))-f_{t,i}(x_{i}(t))\right]\right).\label{eq:th1_2}
\end{align}
Based on the convexity of $f_{t,i}$'s, we have 
\begin{align}
\sum_{t=1}^{T}\left(\frac{1}{n}\sum_{i=1}^{n}f_{t,i}(x_{i}(t))-f_{t}(x^{*})\right)\quad & \leq\nonumber \\
\sum_{t=1}^{T}\left(\frac{1}{n}\sum_{i=1}^{n}\left\langle g_{i}(t),x_{i}(t)-x^{*}\right\rangle \right),\label{eq:th1_3}
\end{align}
where $g_{i}(t)\in\partial f_{t,i}(x_{i}(t))$ is the sub-gradient
of $f_{t,i}$ at $x_{i}(t)$. Thereby, we can express the regret bound
based on \eqref{eq:th1_2}, \eqref{eq:th1_3}, and the $L$-Lipschitz
continuity of $f_{t,i}$'s as, 
\begin{align}
R_{T}(x^{*},x_{i}) & \leq\sum_{t=1}^{T}(\frac{1}{n}\sum_{i=1}^{n}\left\langle g_{i}(t),x_{i}(t)-x^{*}\right\rangle \nonumber \\
 & +\frac{L}{n}\sum_{i=1}^{n}\Vert x_{i}(t)-\phi(t)\Vert+L\Vert x_{i}(t)-\phi(t)\Vert).\label{eq:th1_4}
\end{align}
The first term on the right had side of \eqref{eq:th1_4} can be expanded
as 
\begin{align}
\sum_{t=1}^{T}\left(\frac{1}{n}\sum_{i=1}^{n}\left\langle g_{i}(t),x_{i}(t)-x^{*}\right\rangle \right)\qquad\nonumber \\
=\sum_{t=1}^{T}(\frac{1}{n}\sum_{i=1}^{n}\left\langle g_{i}(t),x_{i}(t)-\phi(t)\right\rangle \nonumber \\
+\frac{1}{n}\sum_{i=1}^{n}\left\langle g_{i}(t),\phi(t)-x^{*}\right\rangle ).\qquad\label{eq:gi expansion}
\end{align}
Now, we need to bound the terms on the right hand side of \eqref{eq:gi expansion}.
The first term is bounded based on the convexity and $L$-Lipschitz
continuity of $f_{t,i}$.%
\footnote{Note that convexity of $f_{t,i}$ implies $\left\langle g_{i}(t),x-y\right\rangle \leq f_{t,i}(x)-f_{t,i}(y)$
. Therefore, based on $L$-Lipschitz continuity of $f_{t,i}$'s, we
have $||g_{i}||_{*}\leq L$ and we can deduce \eqref{eq:<gi,xi-y> bound}.%
} In other words,

\begin{equation}
\left\langle g_{i}(t),x_{i}(t)-\phi(t)\right\rangle \leq L\Vert x_{i}(t)-\phi(t)\Vert.\label{eq:<gi,xi-y> bound}
\end{equation}
Since $x_{i}(t)$ and $\phi(t)$ are the projections of $y_{i}(t)$
and $\bar{y}(t)$ respectively, the Lipschitz continuity of $\Pi_{\chi}^{\psi}(.,\alpha)$
presented in Lemma \ref{lem:x_i-y-1} of the Appendix imposes a bound
on $\Vert x_{i}(t)-\phi(t)\Vert$ as 
\begin{equation}
\Vert x_{i}(t)-\phi(t)\Vert\leq\alpha(t)\Vert\bar{y}(t)-y_{i}(t)\Vert_{*},\label{eq:xi-y bound}
\end{equation}
where $||.||_{*}$ is the dual norm. Therefore, using the bound in
Lemma \ref{lem:y-x*-1} of the Appendix and noting that $\Vert g_{i}(t)\Vert_{*}\leq L$,
we can write \eqref{eq:gi expansion} as 
\begin{align}
\sum_{t=1}^{T}\left(\frac{1}{n}\sum_{i=1}^{n}\left\langle g_{i}(t),x_{i}(t)-x^{*}\right\rangle \right)\quad\nonumber \\
\leq\frac{L}{n}\sum_{t=1}^{T}\sum_{i=1}^{n}\alpha(t)\Vert\bar{y}(t)-y_{i}(t)\Vert_{*}\quad\nonumber \\
+\frac{L^{2}}{2}\sum_{t=2}^{T}\alpha(t-1)+\frac{1}{\alpha(T)}\psi(x^{*}).\label{eq:<gi,xi-x*> bound}
\end{align}
Thus, \eqref{eq:th1_4}, \eqref{eq:xi-y bound}, and \eqref{eq:<gi,xi-x*> bound}
imply that 
\begin{align}
R_{T}(x^{*},x_{i})\leq\frac{L^{2}}{2}\sum_{t=2}^{T}\alpha(t-1)+\frac{1}{\alpha(T)}\psi(x^{\star})\qquad\quad\quad\nonumber \\
+L\sum_{t=1}^{T}\alpha(t)\left(\Vert\bar{y}(t)-y_{i}(t)\Vert_{*}+\frac{2}{n}\sum_{i=1}^{n}\Vert\bar{y}(t)-y_{i}(t)\Vert_{*}\right).\label{eq:th1_5}
\end{align}
On the other hand, Lemma \ref{lem:z-z_i-1} of the Appendix imposes
an upper bound on the last term on the right hand side of \eqref{eq:th1_5}.
Thus, using \eqref{eq:network effedct 3} the regret is further bounded
as 
\begin{align}
R_{T}(x^{*},x_{i}) & \leq\frac{L^{2}}{2}\sum_{t=1}^{T-1}\alpha(t)+\frac{1}{\alpha(T)}\psi(x^{*})\nonumber \\
 & +3L^{2}\left(\frac{n}{1-\gamma}+n\delta\nu+2\left(1-n\right)\right)\sum_{t=1}^{T}\alpha(t).\label{eq:th1_6}
\end{align}
The statement of the theorem now follows from the integral test on
$\alpha(t)=k/\sqrt{t}$ and $\psi(x^{*})\leq R^{2}$. 
\end{IEEEproof}
Theorem \ref{thm:Regret(xi)} indicates the ``good'' performance
of online-DWDA through sub-linear regret and highlights the importance
of the underlying network topology through the parameters $\gamma$
and $\nu$ examined in \S \ref{sec:Convergence-Analysis}. 

Next we present the regret analysis for the (temporal) running average
estimates at each agent exhibiting a similar dependence on the network
connectivity. 
\begin{cor}
Given the sequence\textup{ $\widetilde{x}_{i}(t)$} generated by line
\textup{\ref{eq:x_avg update}} in \textup{Algorithm \ref{alg:DWDA}}
for all $i\in[n]$ with $\psi(x^{*})\leq R^{2}$ and $\alpha(t)=k/\sqrt{t}$,
we have 
\[
R_{T}(x^{*},\widetilde{x}_{i})\leq2\left(\frac{R^{2}}{k}+kL^{2}\left(\frac{6n}{1-\gamma}+6n\delta\nu+1\right)\right)\sqrt{T}.
\]
\end{cor}
\begin{IEEEproof}
Since the cost function $f_{t}(x(t))$ is convex, $f_{t}(\widetilde{x}_{i}(t))\leq\frac{1}{t}\sum_{s=1}^{t}f_{t}(x_{i}(s))$.
Therefore, we have 
\begin{equation}
f_{t}(\widetilde{x}_{i}(t))-f_{t}(x^{*})\leq\frac{1}{t}R_{t}(x^{*},x_{i}).\label{eq:avg regret 2}
\end{equation}
Thus, the running average regret is bounded as 
\begin{equation}
R_{T}(x^{*},\widetilde{x}_{i})\leq\sum_{t=1}^{T}\left(\frac{1}{t}R_{t}(x^{*},x_{i})\right).\label{eq:avg regret vs regret}
\end{equation}
On the other hand, the regret bound \eqref{eq:regret bound} implies
that 
\begin{align}
R_{T}(x^{*},\widetilde{x}_{i}) & \leq\left(\frac{R^{2}}{k}+kL^{2}\left(\frac{6n}{1-\gamma}+6n\delta\nu+1\right)\right)\nonumber \\
 & \;\;\;\,\times\left(\sum_{t=1}^{T}\frac{1}{\sqrt{t}}\right).\label{eq:avg regret 3}
\end{align}
The statement of the corollary now follows from the integral test
on the right hand side of \eqref{eq:avg regret 3}. 
\end{IEEEproof}

\section{Online Distributed Estimation\label{sec:Simulation results}}

\begin{figure}[tp]
\noindent \begin{centering}
\includegraphics[width=1\linewidth]{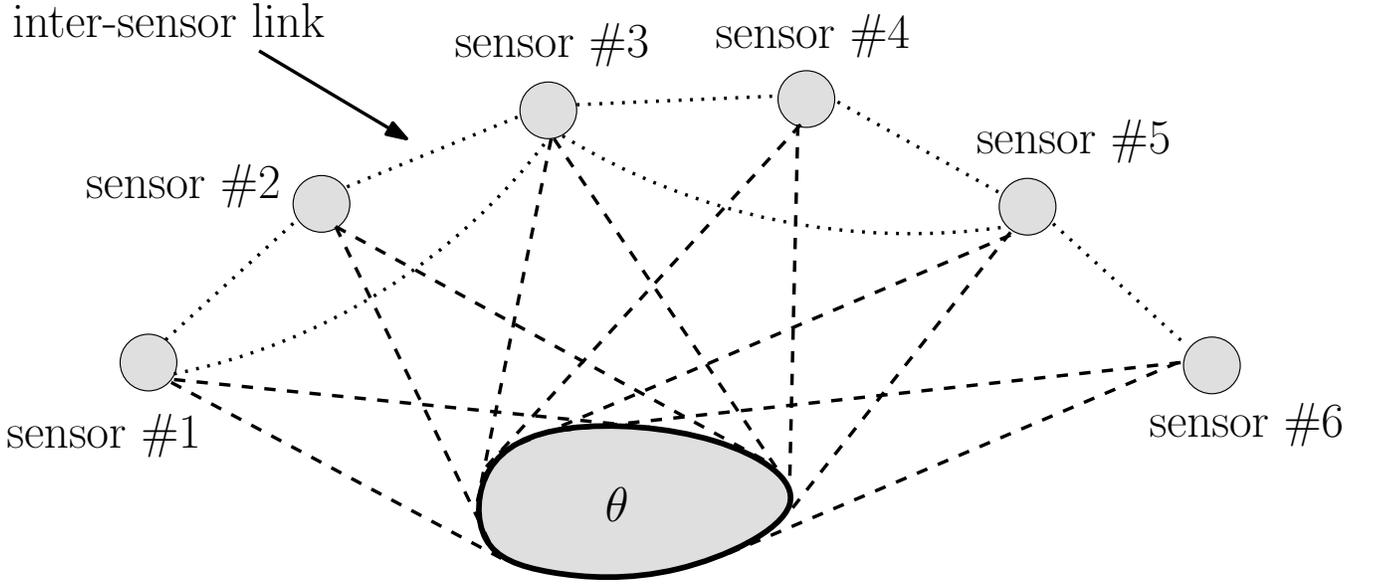} 
\par\end{centering}

\protect\caption{\label{fig:Distributed Network Representation}A graphical representation
of a distributed sensor network.}
\end{figure}

Adopting the least squares point of view, a model for online estimation
over a distributed sensor network is presented in this section. The
distributed sensor network aims to estimate a random vector $\theta\in\Theta=\left\{ \theta\in\mathbb{R}^{d}|\left\Vert \theta\right\Vert _{2}\leq\theta_{\max}\right\} $.
Note that $\Theta$ is a closed convex set containing the origin.
The observation vector $z_{t,i}:\mathbb{R}^{d}\rightarrow\mathbb{R}^{p_{i}}$
represents the $i$th sensor measurement at time $t$ which is uncertain
and time-varying due to the sensor's susceptibility to unknown environmental
factors such as jamming. The sensor is assumed (not necessarily accurately)
to have a linear model of the form $h_{i}(\theta)=H_{i}\theta$, where
$H_{i}\in\mathbb{R}^{p_{i}\times d}$ is the observation matrix of
sensor $i$ and $\left\Vert H_{i}\right\Vert _{1}\leq h_{\max}$ for
all $i$. Consider now the interconnection topology between the sensors
defined via the directed graph $\mathcal{G}=\left(V,E,W\right)$,
where the set of $n$ sensors are represented by $V$. The presence
of an edge $(j,i)\in E$ indicates an information flow from sensor
$j$ to sensor $i$. The set of agents that are communicating with
agent $i$ is defined as the neighborhood set $N(i)=\{j\in V|(j,i)\in E\}$.
Figure \ref{fig:Distributed Network Representation} graphically summarizes
the problem setup. The objective is to find the argument $\hat{\theta}$
that minimizes the cost function 
\begin{eqnarray}
f_{t}(\hat{\theta})=\frac{1}{n}\sum_{i=1}^{n}f_{t,i}(\hat{\theta}) & \;\mbox{subject to}\;\hat{\theta}\in\Theta,\label{eq:network estimation cost}
\end{eqnarray}
where 
\begin{equation}
f_{t,i}(\hat{\theta})=\frac{1}{2}\left\Vert z_{i,t}-H_{i}\hat{\theta}\right\Vert _{2}^{2}\label{eq:sensor i cost}
\end{equation}
is a convex cost function associated with sensor $i\in[n]$. It is
assumed that the value of this local cost at time $t$ is only revealed
to the sensor after $\hat{\theta}(t)$ has been computed, that is,
the local error functions are allowed to change over time in an unpredictable
manner due to modeling errors and uncertainties in the environment.
The (sub)gradient of the local estimation error \eqref{eq:sensor i cost},
\begin{equation}
\partial f_{t,i}(\hat{\theta})=H_{i}^{T}\left(z_{t,i}(\theta)-H_{i}\hat{\theta}\right),\label{eq:sensor i gradient}
\end{equation}
is also assumed to be known to the sensor and its neighbors. We note
that the cumulative cost at time $T$ is defined as $f(\hat{\theta})=\sum_{t=1}^{T}f_{t}(\hat{\theta}).$

\textcolor{black}{In an offline setting, for all $t\in[T]$, each
sensor $i$ has a noisy observation} $z_{t,i}=H_{i}\theta+v_{t,i}$,
where $v_{t,i}$ is generally assumed to be (independent) white noise.
In this case, the centralized time-averaged optimal estimate for \eqref{eq:network estimation cost}
is 
\begin{equation}
\theta^{*}=\frac{1}{T}\sum_{t=1}^{T}\left(\sum_{i=1}^{n}H_{i}^{T}\Sigma_{t,i}^{-1}H_{i}\right)^{-1}\left(\sum_{i=1}^{n}H_{i}^{T}\Sigma_{t,i}^{-1}z_{i,t}\right),\label{eq:best fixed estimation}
\end{equation}
where $\Sigma_{t,i}$ is the covariance of the error observed by sensor
$i$ at time $t$ \cite{Xiao2006}. For the case where $\theta\in\mathbb{R}$,
$\Sigma_{t,i}=I$, and $H_{i}=1$, the optimal estimate is $\theta^{*}=\frac{1}{nT}\sum_{i=1}^{n}\sum_{t=1}^{T}z_{t,i}.$
However, this approach to estimation problems is not suitable in scenarios
where the noise characteristics are unknown. For example when a wireless
sensor network is employed in an unknown and dynamic environment,
the measurement signal can be blocked or degraded due to obstructions
such as walls, furniture, trees, or buildings. This is known as the
shadowing effect and usually modeled as a function of the environment
in which the network is deployed. Another example is jamming of one
or more sensors in the network. When the sensor resolution and noise
characteristics are not known ahead of time, the dynamic weight selection
procedure discussed in \S\ref{sec:Distributed-Adaptive-Weight} can
be employed to eliminate the information from the jammed sensors. 

An online framework is particularly suitable for such estimation problems
without relying on prior assumption or knowledge of the statistical
properties of the data. In the proposed distributed estimation algorithm,
at time step $t$, each sensor $i$ estimates $\hat{\theta}_{i}\in\Theta$
based on the local information available to it and then an ``oracle''
reveals the cost $f_{t}(\hat{\theta}_{i})$.

The bounds presented in Theorem \ref{thm:Regret(xi)} apply after
selecting $\psi(\hat{\theta})=\frac{1}{2}\|\hat{\theta}\|_{2}^{2}$
and the parameter $\alpha(t)$ accordingly. In order to find the constants
$R$ and $L$ featured in the result, we note that for $\hat{\theta}\in\Theta$,
$\psi(\hat{\theta})\leq\frac{1}{2}\theta_{\max}^{2}$, and thus $R\leq\frac{1}{\sqrt{2}}\theta_{\max}$.
In this example, we assume that the observation for agent $i$ at
time $t$ is of the form $z_{t,i}=a_{t}\theta+b_{t}$ for some $a\in\left(0,a_{\max}\right)$
and $b\in\left(-b_{\max},b_{\max}\right)$. Therefore, 
\[
\sup_{\theta\in\chi}\left\Vert z_{t,i}(\theta)\right\Vert _{2}\leq a_{\max}\theta_{\max}+b_{\max}.
\]
Further, the function $f_{t,i}$ is Lipschitz as it is convex on a
compact domain and the Lipschitz constant can be found by observing
that 
\begin{eqnarray*}
 &  & \left|f_{k,i}(\hat{\theta})-f_{k,i}(\phi)\right|\\
 &  & \leq\frac{1}{2}\left|\left(\hat{\theta}-\phi\right)^{T}H_{i}^{T}H_{i}\left(\hat{\theta}-\phi\right)\right|+\left|z^{T}H_{i}\left(\hat{\theta}-\phi\right)\right|\\
 &  & \leq\left(\frac{1}{2}\left\Vert H_{i}\right\Vert _{F}^{2}\left\Vert \hat{\theta}-\phi\right\Vert _{2}+\left\Vert z_{t,i}\right\Vert _{2}\left\Vert H_{i}\right\Vert _{F}\right)\left\Vert \hat{\theta}-\phi\right\Vert _{2}
\end{eqnarray*}
and thus $L=\left(\frac{1}{2}\theta_{\max}h_{\max}+a_{\max}\theta_{\max}+b_{\max}\right)h_{\max}$.
Hence $R_{T}(\theta^{*},\hat{\theta}_{i})/T\rightarrow0$ and the
algorithm performs as well as best fixed estimate $\theta^{*}$ in
hindsight \eqref{eq:best fixed estimation} ``on average''. For
the case where $\theta_{t}=\theta_{t+1}$ for $t=1,2,\dots,T,$ $\theta^{*}$
is the optimal estimate.

The online-DWDA and DOA algorithms have been implemented on the described
distributed sensor setup for $n=100$ sensors. The objective is to
estimate a scalar $\theta\in\left(-\frac{1}{2},\frac{1}{2}\right)$
with a fixed $H_{i}\in\left(0,\frac{1}{4}\right)$ for each agent;
hence $\sup_{i}\left|H_{i}\right|=\frac{1}{4}.$ In this example,
we have assumed $a\in\left(0,1\right)$, $b\in\left(-\frac{1}{4},\frac{1}{4}\right)$,
$\beta=0.9$, and $k=\frac{1}{4}.$ Thus, $d=1,$ $\Theta=\left(-\frac{1}{2},\frac{1}{2}\right)$,
$h_{\max}=\frac{1}{4}$ , $\theta_{\max}=\frac{1}{2}$, $R=\frac{1}{2\sqrt{2}}$,
and $L=\frac{13}{64}$.

\begin{figure}[tp]
\noindent \begin{centering}
\includegraphics[width=1\linewidth]{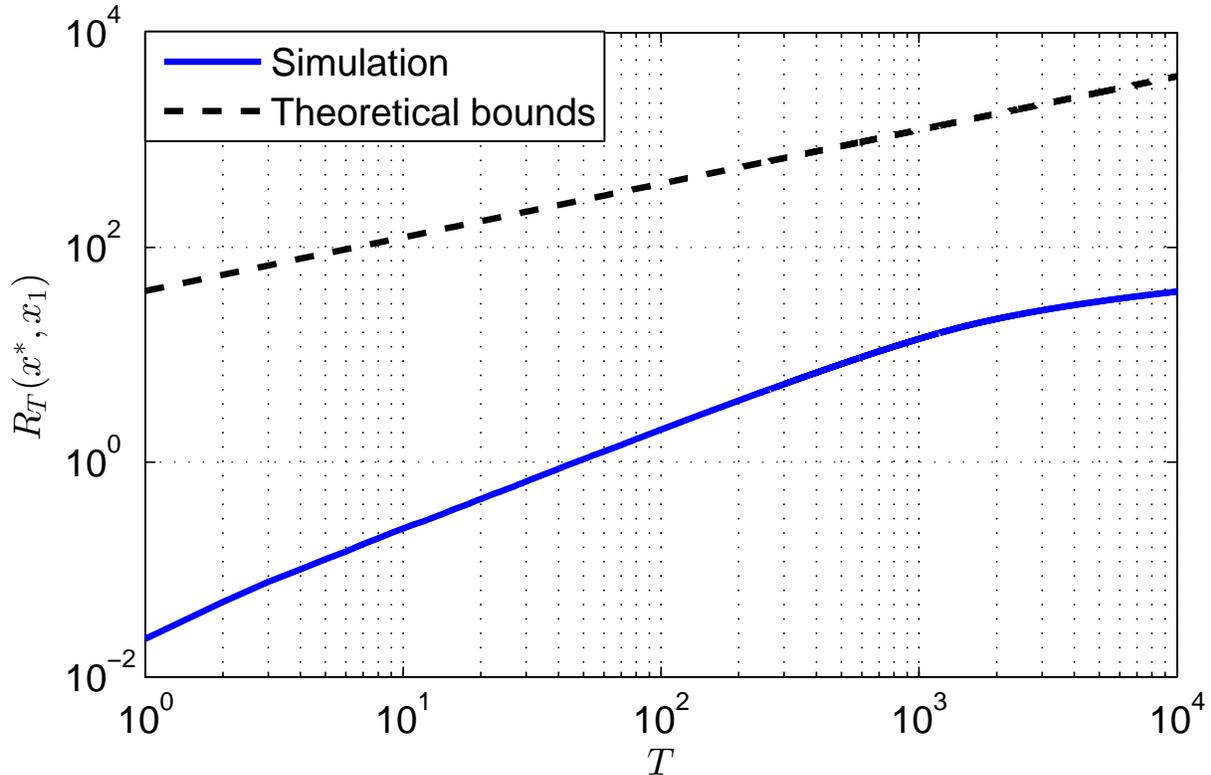} 
\par\end{centering}

\protect\caption{\label{fig:BoundAccuracy}Accuracy of the bounds in \eqref{eq:regret bound}
where $\mathcal{G}$ is a 100 node random directed graph with edge
probability $p=0.08$, $\nu=5$, and $\gamma=0.2034$.}
\end{figure}

\begin{figure}[tp]
\noindent \begin{centering}
\includegraphics[width=1\linewidth]{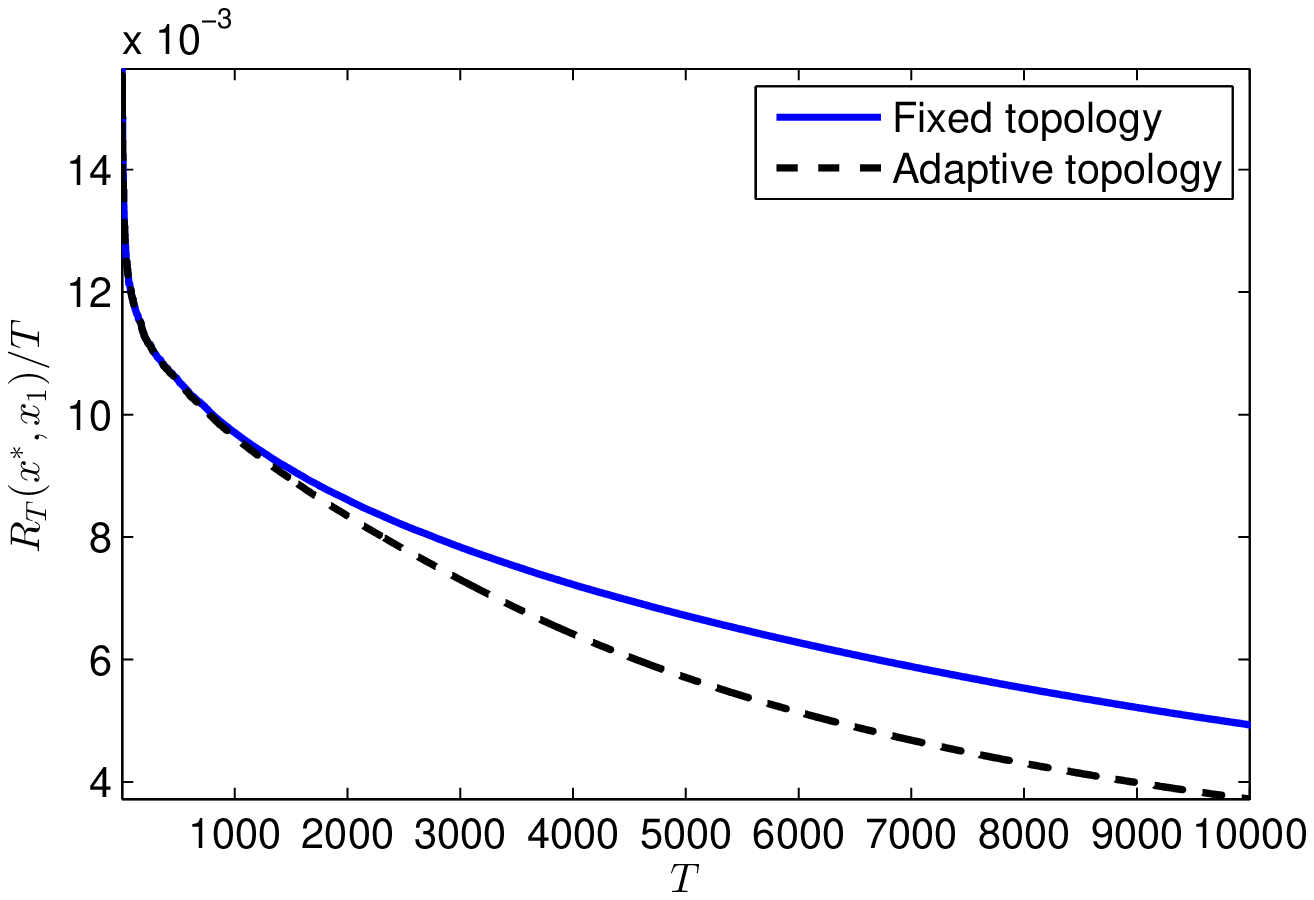} 
\par\end{centering}

\protect\caption{\label{fig:Sensor jamming}Regret performance over fixed and adaptive
network topologies where $\mathcal{G}$ is a 100 node random $4$-regular
graph and $25$ sensors are assumed to have been jammed. For the jammed
sensors, $b_{t}=b_{\max}$ and $a_{t}=H_{i}$.}
\end{figure}

The online-DWDA and DOA algorithms were also applied to random sensor
network with edge probability $p=0.08$. Figure \ref{fig:BoundAccuracy}
shows a qualitative agreement of the theoretical regret bound \eqref{eq:regret bound}
and simulation results, indicating that $R_{T}(x^{*},x_{1})=O(\sqrt{T})$.
The improved performance of the adaptive network topology has been
emphasized in Figure \ref{fig:Sensor jamming} in the context of a
jamming scenario, where a number of sensors in the random regular
network are assumed to have been jammed. This figure also demonstrates
that the adaptive sensor network has a better regret performance as
compared with the fixed topology sensor network.

\begin{figure}[tp]
\noindent \begin{centering}
\includegraphics[width=1\linewidth]{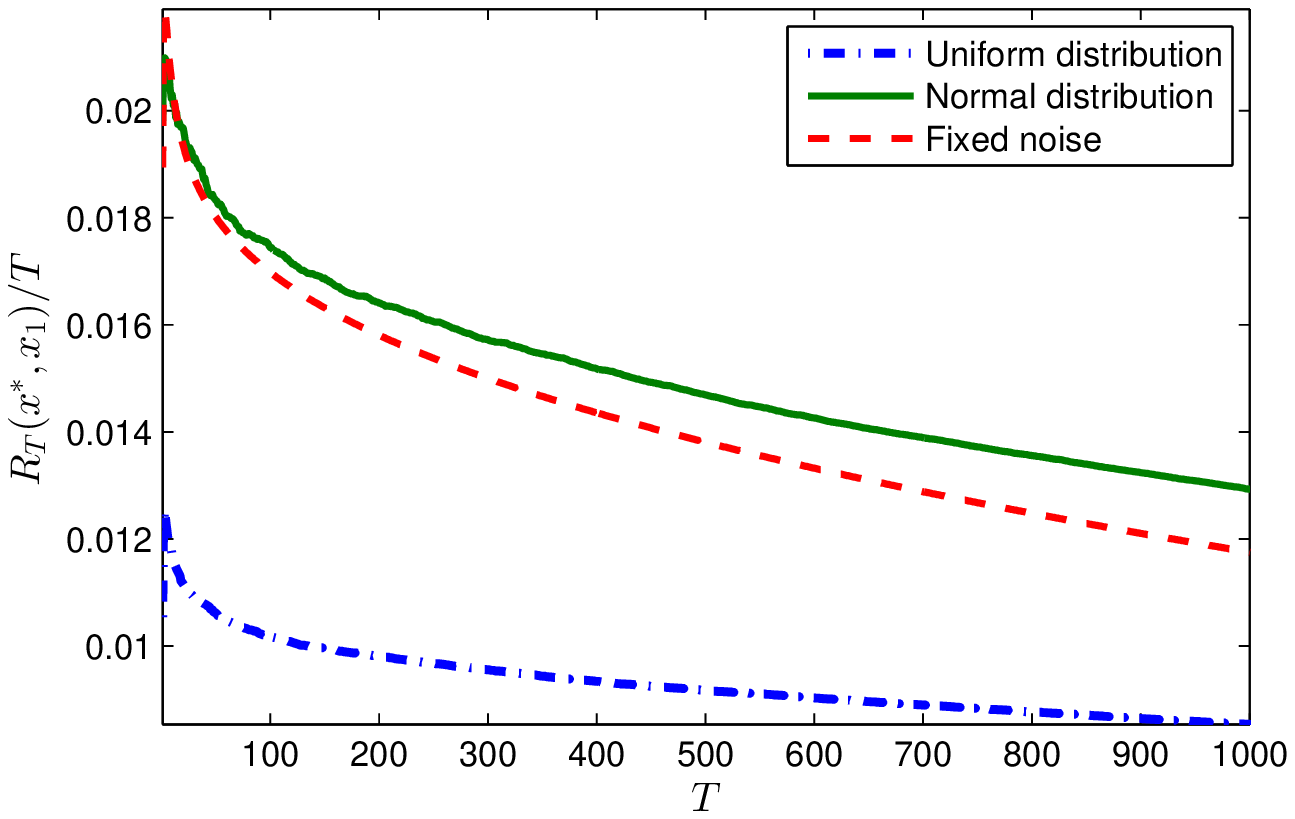} 
\par\end{centering}

\protect\caption{\label{fig:Regret vs observation noise}Regret performance for three
different observation noise characteristics, where $\mathcal{G}$
is a 100 node random $4$-regular graph. The noise signals have been
generated from distributions with mean $-b_{\max}$ and standard deviation
$b_{\max}$.}
\end{figure}

In addition, the performance of the proposed adaptive online distributed
estimation in the presence of various noise types is presented in
Figure \ref{fig:Regret vs observation noise}. These simulation results
indicate that $R_{T}(\theta^{*},\hat{\theta}_{1})=O(\sqrt{T})$ for
all noise types considered without a prior assumption on the noise
characteristics. 

\begin{figure}[tp]
\noindent \begin{centering}
\includegraphics[width=1\linewidth]{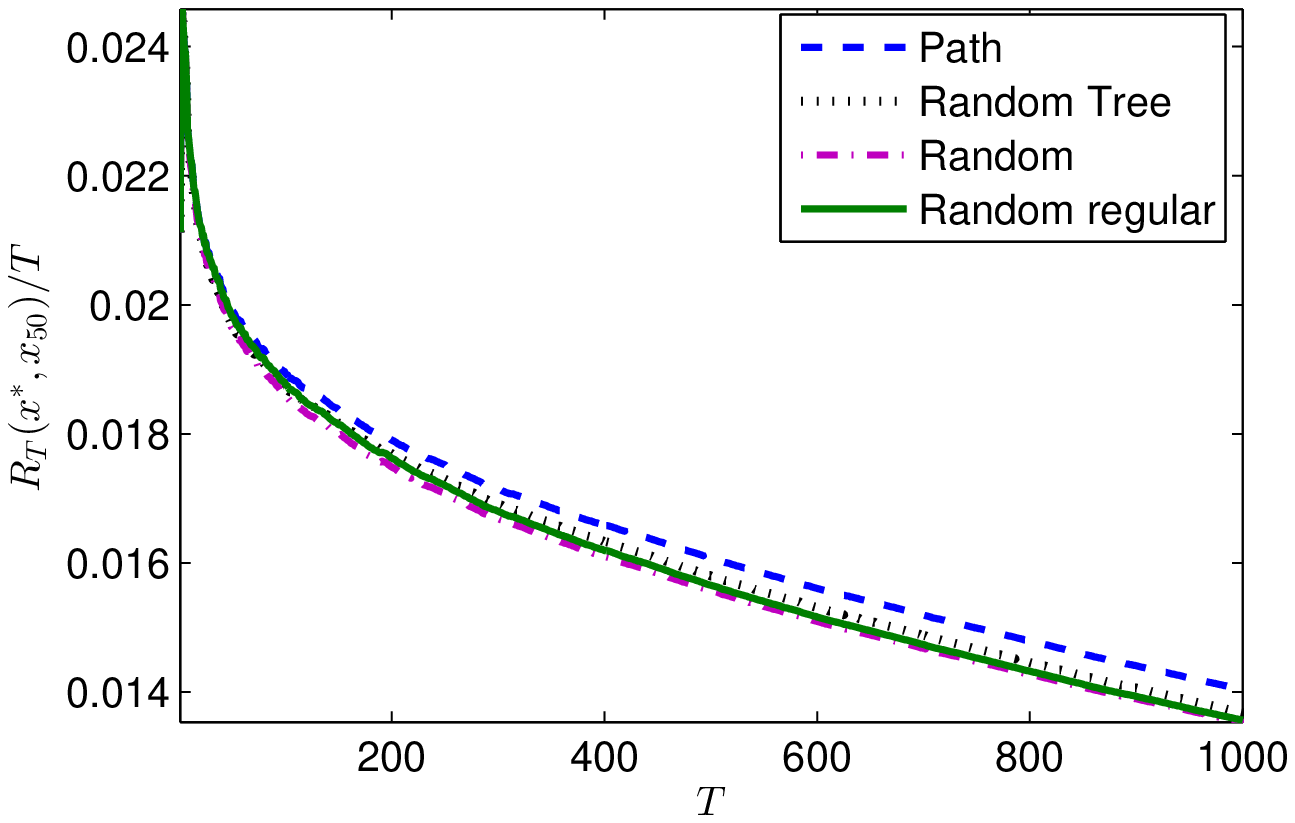} 
\par\end{centering}

\protect\caption{\label{fig:GraphType} The performance of the online distributed estimation
algorithm on four different 100 node graphs with $\gamma=\left\{ 0.8999,\;0.8998,\;0.7939,\;0.4110\right\} $,
for the path, random tree, random $k$-regular with $k=4$, and random
graph with edge probability $p=0.08$, respectively, presented in
increasing order of performance.}
\end{figure}

Furthermore, the role of network connectivity in the performance of
the algorithm has been emphasized in Figure \ref{fig:GraphType} for
various classes of network topologies, directly correlated to the
network connectivity measure $\gamma$. This result can be applied
to designing sensor network topologies that operate in highly uncertain
environments. Suitable metrics for such a topology design procedure
include $\lambda_{2}(P(\mathcal{G}^{0}))$ that predictably scales
with $n,$ such as random regular graphs and expander graphs \cite{Bollobas1998}.

\section{Conclusion\label{sec:Conclusion}}

This paper studies the problem of decentralized optimization on dynamic
networks operating in an uncertain environment. An algorithm has been
presented that evolves distributively using only local information
available to the agents in the network. Our analysis provided a convergence
rate of $O(1/\sqrt{T})$ and a sub-linear regret of $O(\sqrt{T})$
in the online setting. In addition, the convergence analysis of the
distributed optimization algorithm highlighted the role of two measures
of network connectivity. 

A distributed dynamic weight selection procedure has also been proposed
that on average, performs as well as the best strategy for information
diffusion in hindsight. It was demonstrated that this approach improves
the convergence rate by mitigating the network effects.

\begin{comment}
In order to justify the suitability of the online setting for sensor
networks, the proposed algorithm was applied to a distributed sensor
estimation problem. The estimates were acquired in real-time and coupled
with sensors' susceptibility to unknown errors and jamming. The simulation
results indicate that the sensor network can provide an estimate that
on the average performs as well as the best case fixed solution in
hindsight. 
In addition, we explored the proposed online distributed estimation
algorithm for various classes of sensor networks and highlighted the
role of network connectivity on the network-level regret. 
\end{comment}

This work can be applied in the context of a range of applications
such as mobile sensor networks where the network is susceptible to
unknown errors, jamming, link failure, and a varying network topology.
Moreover, this work can be extended in several directions. One such
extension, which is the subject of our future work, involves examining
online distributed filtering. More generally, the online approach
can be adopted for a host of network dynamic systems that operate
in unstructured environments, requiring that a learning algorithm
is embedded in the network-level decision-making process.

\section{Appendix\label{sec:Appendix}}

We note that Lemmas \ref{lem:x_i-y-1} and \ref{lem:y-x*-1} have
been shown by Duchi \textit{et al.}, \cite{Duchi2012} and are presented
here for reference. 
\begin{lem}
\cite{Duchi2012}\label{lem:x_i-y-1} For any $u,v\in\mathbb{R}^{m}$,
and under the conditions stated for proximal function $\psi$ and
step size $\alpha(t),$ we have $\Vert\Pi_{\chi}^{\psi}(u,\alpha)-\Pi_{\chi}^{\psi}(v,\alpha)\Vert\leq\alpha\Vert u-v\Vert_{*}$. 
\end{lem}

\begin{lem}
\cite{Duchi2012}\label{lem:y-x*-1} For any positive and non-increasing
sequence $\alpha(t)$ and $x^{*}\in\chi$, 
\[
\sum_{t=1}^{T}\langle\overline{g}(t),\phi(t)-x^{*}(t)\rangle\leq\frac{1}{2}\sum_{t=1}^{T}\alpha(t-1)\Vert\overline{g}(t)\Vert_{\star}^{2}+\frac{1}{\alpha(T)}\psi(x^{*}),
\]
where the sequence $\phi(t)$ is generated by \eqref{eq:y update}. 
\end{lem}
The following result presents a bound on $\Vert\bar{y}(t)-y_{i}(t)\Vert_{*}$
proportional to the error incurred by the decentralized update in
Algorithm \ref{alg:DWDA}. 
\begin{lem}
\label{lem:z-z_i-1}For sequences $y_{i}(t)$ and $\bar{y}(t)$ generated
by line \textup{\ref{eq:z_i update}} of Algorithm\textup{ \ref{alg:DWDA}}
and \eqref{eq:z_bar update}\textup{,} respectively, we have, 
\[
\Vert\bar{y}(t)-y_{i}(t)\Vert_{*}\leq L\sum_{k=0}^{t-2}\sum_{j=1}^{n}\left|P_{ij}^{(t-1,k+1)}-\pi_{j}\right|+2L,
\]
 for all $i\in[n]$.\end{lem}
\begin{IEEEproof}
Reformulating the update in line \ref{eq:z_i update} of Algorithm
\ref{alg:DWDA} for all $i\in[n]$, by induction through $s$ steps
we have, 
\begin{align}
y_{i}(t) & =\sum_{j=1}^{n}P_{ij}^{(t-1,t-s)}y_{j}(t-s)+\sum_{k=t-s}^{t-2}\sum_{j=1}^{n}P_{ij}^{(t-1,k+1)}g_{j}(k)\nonumber \\
 & \;\;+g_{i}(t-1).\label{eq:z_i evol-1}
\end{align}

Since $\bar{y}(t)$ evolves as in \eqref{eq:z_bar update}, by setting
$s=t$ in \eqref{eq:z_i evol-1} and assuming $y_{i}(0)=0$, we get,
\begin{align}
\bar{y}(t)-y_{i}(t) & =\sum_{k=0}^{t-2}\left(\sum_{j=1}^{n}\left(\pi_{j}-P_{ij}^{(t-1,k+1)}\right)g_{j}(k)\right)\nonumber \\
 & \;\;\;+\bar{g}(t-1)-g_{i}(t-1).\label{eq:z-z_i evol-1}
\end{align}

Thus, the dual norm of \eqref{eq:z-z_i evol-1} is bounded as

\begin{align}
\Vert\bar{y}(t)-y_{i}(t)\Vert_{*} & \leq\Vert\sum_{k=0}^{t-2}\left(\sum_{j=1}^{n}\left(\pi_{j}-P_{ij}^{(t-1,k+1)}\right)g_{j}(k)\right)\Vert_{*}\nonumber \\
 & +\Vert\bar{g}(t-1)-g_{i}(t-1)\Vert_{*},\label{eq: dual (z-z_i)-1}
\end{align}
and the right hand side of \eqref{eq: dual (z-z_i)-1} can be bounded
by

\begin{align}
\Vert\bar{y}(t)-y_{i}(t)\Vert_{*} & \leq\sum_{k=0}^{t-2}\sum_{j=1}^{n}\left|P_{ij}^{(t-1,k+1)}-\pi_{j}\right|\Vert g_{j}(k)\Vert_{*}\nonumber \\
 & +\Vert\bar{g}(t-1)-g_{i}(t-1)\Vert_{*}.\label{eq:dual (z-z_i) 2-1}
\end{align}
Since $\Vert g_{i}(t)\Vert_{\star}\leq L$,
\begin{equation}
\Vert\bar{y}(t)-y_{i}(t)\Vert_{*}\leq L\sum_{k=0}^{t-2}\sum_{j=1}^{n}\left|P_{ij}^{(t-1,k+1)}-\pi_{j}\right|+2L.\label{eq:dual (z-z_i) 3-1}
\end{equation}

\end{IEEEproof}
The following proposition provides an upper bound on $\nu$ in the
convergence rate \eqref{eq:convergence rate} over fixed topology
networks.
\begin{prop}
\label{prop: shortest path}Consider a set $\mathcal{P}$ of stochastic
matrices with positive diagonal elements, representing an arbitrarily
strongly connected topologies over $n$ nodes, i.e., $P^{k}\in\mathcal{P}$
for all positive integers $k$. Suppose that any two matrices $P^{k_{1}}$
and $P^{k_{2}}$ are of the same type.%
\footnote{The matrices $A$ and $B$ are of the same type if they have zero
elements and positives elements in the same place.%
} Then, there exists an integer $\nu$, 
\begin{equation}
1\leq\nu\leq\min_{i\in[n]}\max_{j\in[n]}\mbox{dist}(j,i),\label{eq:v bound}
\end{equation}
for which if the sequence $Q=P^{(m+\nu-2,m)}$ of matrices in $\mathcal{P}$
is not scrambling, $P^{m+\nu-1}Q$ is scrambling. \end{prop}
\begin{IEEEproof}
Let $Q_{1}=P^{m}$ and $Q_{2}=P^{m+1}P^{m}.$ Thus every entry of
$Q_{2}$ is represented as 
\[
\left[Q_{2}\right]_{ij}=\sum_{k=1}^{n}\left[P^{m+1}\right]_{ik}\left[P^{m}\right]_{kj}.
\]
Since $\left[P^{m}\right]_{ii}>0$ for all $i\in[n]$ and integer
$m\geq1$, the entry $\left[Q_{2}\right]_{ij}$ is positive if $\left(j,i\right)\in E,$
$\left(i,j\right)\in E$, or if there exists a node $k\in[n]$ in
the directed path from node $j$ to node $i$ with $\mbox{dist}(j,i)=2.$
Thus, the corresponding zero entry of $Q_{1}$ that has one of the
aforementioned properties will be positive in $Q_{2}.$ By induction,
it follows that the entry of $\left[Q_{\nu}\right]_{ij}$ will be
positive if $\left(j,i\right)\in E$, $\left(i,j\right)\in E$, or
if there exists a node $k\in[n]$ in the directed path from node $j$
to node $i$ with $\mbox{dist}(j,i)=\nu_{i}$. Therefore, for each
row $i$ of $Q_{\nu}$, all entries will be positive when 
\[
\nu_{i}=\max_{j\in\left[n\right]}\mbox{dist}(j,i).
\]
Note that every element of any row of the sequence $Q=P^{(m+\nu-2,m)}$
of matrices in $\mathcal{P}$ is positive, the matrix $Q$ is scrambling
and $\nu$ satisfies the bound \eqref{eq:v bound}. 
\end{IEEEproof}
A similar observation for the adjacency matrix of $\mathcal{G}$ can
be found in the algebraic graph theory literature such as \cite{Godsil2001}. 

\bibliographystyle{IEEEtran}
\bibliography{generalPapers,regretPapers,NetworkTopology,DSSLOnlineRegretCDC2013,DSSLOnlineRegretADMM,DSSLOnlineRegret,cdc2013distopt,Distributedestimation}

\end{document}